\documentclass{amsart}

\usepackage{amssymb}

\newtheorem{thm}{Theorem}

\newtheorem{question}[thm]{Question}

\newtheorem{prop}[thm]{Proposition}

\theoremstyle{definition}
\newtheorem{defn}[thm]{Definition}

\newtheorem{say}[thm]{}

\newtheorem{rem}[thm]{Remark}          

\newtheorem{ack}{Acknowledgments}

\newtheorem{defn-thm}[thm]{Definition--Theorem}  
\newtheorem{defn-lem}[thm]{Definition--Lemma}  

\theoremstyle{remark}

\renewcommand{\c}[0]{{\mathbb C}}

\newcommand{\z}[0]{{\mathbb Z}}

\newcommand{\p}[0]{{\mathbb P}}
\newcommand{\f}[0]{{\mathbb F}}
\newcommand{\q}[0]{{\mathbb Q}}
\newcommand{\map}[0]{\dasharrow}
\newcommand{\qtq}[1]{\quad\mbox{#1}\quad}
\newcommand{\spec}[0]{\operatorname{Spec}}
\newcommand{\pic}[0]{\operatorname{Pic}}
\newcommand{\gal}[0]{\operatorname{Gal}}

\newcommand{\aut}[0]{\operatorname{Aut}}

\newcommand{\chow}[0]{\operatorname{Chow}}
\newcommand{\chr}[0]{\operatorname{char}}

\newcommand{\hilb}[0]{\operatorname{Hilb}}

\newcommand{\bir}[0]{\operatorname{Bir}}    
\newcommand{\graph}[0]{\operatorname{graph}}    
\newcommand{\isom}[0]{\operatorname{Isom}}

\begin{document}
\bibliographystyle{amsalpha}

\title[Birational rigidity of Fano varieties]
{Birational rigidity of Fano varieties \\
and field extensions}
\author{J\'anos Koll\'ar}

\maketitle


The modern study of the birational properties of Fano varieties
started with the works of Iskovskikh; see the surveys
\cite{isk, che1} and the many references there.
A key concept that emerged in this area is birational rigidity.

 Let $X$ be a Fano variety with $\q$-factorial, terminal singularities
and Picard number 1. Roughly speaking, $X$ is called
birationally rigid if $X$ can not be written
in terms of Fano varieties in any other way. 
The precise definition is the following.

\begin{defn}\label{first.defn}
A {\it Mori fiber space} is a
projective morphism $f:X\to Y$ such that
\begin{enumerate}
\item $X$ is $\q$-factorial with terminal singularities,
\item the relative Picard number of $X/Y$ is 1, and
\item $-K_X$ is $f$-ample.
\end{enumerate}

 Let $X$ be a Fano variety with $\q$-factorial, terminal singularities
and Picard number 1 defined over a field $k$. $X$ is called
{\it birationally rigid}  if $X$ is not birational to the
total space of any Mori fiber space, save the trivial
one $f:X\to \spec k$.

For our purposes, it is better to separate this condition
into 2 parts:
\begin{enumerate}\setcounter{enumi}{3}
\item Let $X'$ be a  Fano variety with $\q$-factorial, terminal singularities
and Picard number 1 over $k$ that is birational to $X$. Then
$X'$ is isomorphic to $X$.
\item The only 
rational map $g:X\map Y$ with rationally connected general fibers 
is the constant map $X\map (\mbox{point})$.
\end{enumerate}

In characteristic 0, \cite{bchm}
implies that a map $g:X\map Y$ as in (5)
leads to a Mori fiber space $g':X'\to Y'$
with $X'$ birational to $X$ but possibly $Y'$ not birational to $Y$.
Thus in characteristic 0 the two versions are equivalent
but in positive characteristic this is not known.
\end{defn}

The aim of this note is to settle some foundational questions
about the behavior of birational rigidity in
extensions of algebraically closed fields.
In all cases when the birational rigidity of a variety has been proved,
the untwisting of birational maps 
is done by clear  geometric  constructions
that  are independent of the  algebraically closed field
of definition. Our first result  shows that, 
over  algebraically closed fields, any untwisting behaves similarly.

\begin{thm}\label{bir.thm}
 Let $k$ be an algebraically closed field
and $K\supset k$ an  algebraically closed overfield.
Let $X_k$ be  a
Fano variety  with  Picard number 1 defined over $k$.
Assume that 
\begin{enumerate}
\item either $\chr k=0$ and $X_k$ has $\q$-factorial, terminal singularities
\item or that  $X_k$ is smooth.
\end{enumerate}
 Then  $X_k$ is birationally rigid iff $X_K$ is birationally rigid.
\end{thm}

\begin{rem} I don't know how $\q$-factoriality behaves in
flat families in positive characteristic.
This is the main reason why smoothness is assumed.
For instance, a cone over an elliptic curve is
$\q$-factorial over $\bar{\f}_p$, but not
$\q$-factorial over  any other 
algebraically closed field. This example has
log canonical singularities, but I do not know
how to exclude this phenomenon for Fano varieties
with terminal singularities.
(There is also the slight problem that the definition
of terminal may not be clear in positive characteristic.)
\end{rem}

There is a very interesting problem
related to Theorem \ref{bir.thm}. 

\begin{question} \label{quest}
Let $X$ be a Fano variety over a field $k$
such that $X$ is  birationally rigid over the
algebraic closure $\bar k$. Is $X$ birationally rigid over $k$?
\end{question}

In the terminology of Cheltsov, this asks if
the notions of birational rigidity and
universal birational rigidity coincide or not.

I think that this is very unlikely but I do not have a counter example.
If $G$ denotes the Galois group $\gal(\bar k/k)$, then the
$k$-forms of $X_{\bar k}$ are
classified by $H^1(G,\aut(X_{\bar k}))$ and two such forms are
birational if they have the same image in
$H^1(G,\bir(X_{\bar k}))$. For an arbitrary variety
the map $H^1(G,\aut(X_{\bar k}))\to H^1(G,\bir(X_{\bar k}))$
is not injective. It is quite interesting that
for many birationally rigid varieties, the group
$\bir(X_{\bar k})$ is a split extension
$$
1\to \Gamma\to \bir(X_{\bar k})\to \aut(X_{\bar k})\to 1
$$
for a subgroup $\Gamma$ generated by the
``obvious'' birational self-maps. For all
such examples, $H^1(G,\aut(X_{\bar k}))\to H^1(G,\bir(X_{\bar k}))$
is an injection.
See \cite{mella, che2, shr} for several relevant examples.

Note also that the  question (\ref{quest}) is not
equivalent to the above Galois cohomology problem.
Given $X_k$, in birational geometry we are also interested in birational
equivalences $X_k\sim X'_k$
where $X'_k$ has Picard number 1 over $k$ but
higher Picard number  over $\bar k$.
\medskip

In order to study birational maps of products  of Fano varieties,
\cite{che2} introduced a variant of question (\ref{quest}).
(In the terminology of \cite{che2},
one asks for  varieties for which
$\bir(X)$ ``universally  untwists maximal centers.'')
Our second result characterizes such varieties.
For this we need to define
the dimension of $\bir(X)$.

\begin{defn} Let $X$ be a projective variety
over an algebraically closed field $k$. 
A birational map $\phi:X\map X$ can be identified with the
closure of its graph  $\graph(\phi)\subset X\times X$.
This construction realizes $\bir(X)$ as an open subscheme
of $\hilb( X\times X)$ or of $\chow( X\times X)$.
Let $\graph:\bir(X)\to \hilb( X\times X)$ denote this injection.
In general, $\graph\bigl(\bir(X)\bigr)$  is an at most
 countable union of finite type subschemes.

We can now define the {\it dimension} of  $\bir(X)$ as the
supremum of the dimensions of all irreducible subsets of
$\graph\bigl(\bir(X)\bigr)$.

This representation, however,  is not particularly unique.
Let us call a subset $Z\subset \bir(X)$ {\it constructible}
if $\graph(Z)$  is constructible as a subset of $\hilb( X\times X)$
(that is, 
a finite union of locally closed subvarieties). 
The notion of a constructible subset is independent of 
the birational model of $X$ and the constructible structure
is compatible with the group multiplication in $\bir(X)$.
The dimension of a constructible subset is also well defined.

We can  also define the  dimension of  $\bir(X)$ as the
supremum of the dimensions of all constructible subsets of
$\bir(X)$.
\end{defn}

\begin{thm}\label{thm}  Let $k$ be an algebraically closed field
of characteristic 0.
Let $X$ be a 
Fano variety  with $\q$-factorial, terminal singularities
and Picard number 1 defined over $k$.
The following are equivalent:
\begin{enumerate}
\item $X$ is  birationally rigid and $\dim \bir(X)=0$.
\item $X$ is  birationally rigid and 
$\bir\bigl(X_K\bigr)=\bir X$ for every 
overfield $K\supset k$.
\item $X_K$ is  birationally rigid and 
$\bir\bigl(X_K\bigr)=\bir X$ for every 
overfield $K\supset k$.
\end{enumerate}
If $k$ is uncountable, these are further equivalent to
\begin{enumerate}\setcounter{enumi}{3}
\item $X$ is  birationally rigid and $\bir(X)$ is countable.
\end{enumerate}
\end{thm}

\begin{say}[Comments on  pliability]
The works of Corti suggest that instead of birational rigidity,
one should consider varieties which are birational to 
only finitely many different Mori fiber spaces
up to square equivalence
\cite{co-re, co-me}. These are called {\it  varieties with
finite pliability}.

The proof of Theorem \ref{bir.thm} also shows
that this notion is independent of 
the choice of an algebraically closed field of defintion.

The proof of Theorem \ref{thm} implies that if
$X$ has finite pliability  and $\bir(X)$ is zero dimensional
then any Mori fiber space birational to
$X\times U$ is of the form $X'\times U$ where
$X'$ is a Mori fiber space birational to $X$.
\end{say}

\begin{say}[Proof of Theorem \ref{bir.thm}]\label{pf.of.bir.thm}

Assume first that $X_K$ is birationally rigid.
Let $X'_k\to Y_k$ be a Mori fiber space whose total
space is $k$-birational to $X_k$.
Then  $X'_K\to Y_K$ is a Mori fiber space whose total
space is $K$-birational to $X_K$, thus
$X'_K$ is $K$-isomorphic to $X_K$.

$\isom(X_k, X'_k)$ is a $k$-variety which has a $K$-point.
Thus it also has a $k$-point and so
$X'_k$ is $k$-isomorphic to $X_k$. Thus
$X_k$ is birationally rigid.

The converse  follows from the next, more general, result.
\end{say}

\begin{prop}\label{bir.prop}
 Let $k$ be an algebraically closed field,
$U$ a $k$-variety and $X_U\to U$ a flat family of
 Fano varieties  with terminal singularities.
If $\chr k\neq 0$ then assume in addition that
 $X_U\to U$ is smooth.
Assume that the set 
$$
R(U):=\{u\in U(k): X_u 
\mbox{ is $\q$-factorial,  birationally rigid and $\rho(X_u)=1$}\}
$$
is Zariski dense.
 Then  the geometric generic fiber  $X_K$ of   $X_U\to U$ 
is birationally rigid.
\end{prop}

Proof. It is easy to see that  the 
Picard number of $X_K$ is also 1.
Let $\pi_K:X'_K\to S_K$ be a Mori fiber space
and $\phi_K:X_K\map X'_K$ a birational map.

Possibly after replacing $U$ by a generically finite ramified cover,
we may assume that the above varieties and maps are defined over $U$.
Thus we have
$$
X_U\stackrel{\phi_U}{\map} X'_U \stackrel{\pi_U}{\to} S_U.
$$
If $S_U\to U$ has positive dimensional fibers,
then for almost all $u\in R(U)$, the fiber $X_u$ is
birational to a nontrivial Mori fiber space; a contradiction.
Thus we may assume that $S_U=U$
and the 
Picard number of $X'_K$ is also 1.

Assume for the moment that
 the Picard number of $X'_u$ is also 1 for general $u\in R(U)$.
Then $X_u\cong X'_u$ since $X_u$ is birationally rigid.

Let us now consider the scheme
$\isom\bigl(X_U, X'_U\bigr)$ parametrizing
isomorphisms between the fibers of
$X_U\to U$ and $X'_U\to U$. Since all isomorphisms
preserve the ample anti-canonical class,
$\isom\bigl(X_U, X'_U\bigr)$ is a scheme of finite type
over $U$ (cf.\ \cite[I.1.10]{rc-book}).
 By assumption $\isom\bigl(X_U, X'_U\bigr)\to U$
has nonempty fibers over the dense subset $R(U)\subset U$.
Therefore  $\isom\bigl(X_U, X'_U\bigr)$
dominates $U$ and so the 
geometric generic fiber  $X_K$ of   $X_U\to U$ 
is isomorphic to the 
geometric generic fiber  $X^*_K$ of   $X^*_U\to U$,
as required.

We have seen that 
the Picard number of the geometric generic fiber of
$X'_U\to U$  also 1. In characteristic 0,
this implies that the Picard number of
$X'_u$ is also 1 for general $u\in U(k)$.
(If $X'_U\to U$ is smooth and $k=\c$, then $\pic(X'_u)=H^2(X'_u,\z)$
implies this. The general singular case is treated in 
 \cite[12.1.7]{ko-mo}.) 
In positive characteristic, the topological arguments of
\cite[12.1.7]{ko-mo} do not apply and I do not know
if the rank of the Picard group is a constructible function
for families of Fano varieties.

In our case, the following auxiliary argument does
the trick.

Let $Z\to X_U\times_UX'_U$ be the normalization of the closure of the graph
of $\phi_U$ with projections $p:Z\to X_U$ and
$p':Z\to X'_U$. Let $E\subset Z$  (resp.\ $E'\subset Z$)
be the exceptional divisors of  $p$ (resp.\ $p'$).
Then
$$
\begin{array}{rcl}
\rho(Z_K)&=&\rho(X_K)+\#\{\mbox{irreducible components of $E_K$}\},\qtq{and}\\
&=&\rho(X'_K)+\#\{\mbox{irreducible components of $E'_K$}\}.
\end{array}
$$
Since $\rho(X_K)=\rho(X'_K)=1$, we obtain that 
 $E_K$ and $E'_K$ have the same
number of irreducible components.
Similarly, for general  $u\in U(k)$, $Z_u$ is the graph of
a birational map from $X_u$ to $X'_u$ and
$$
\begin{array}{rcl}
\rho(X_u)-\rho(X'_u)&=&\#\{\mbox{irred.\ comps.\ of $E_u$}\}
- \#\{\mbox{irred.\ comps.\ of $E'_u$}\}\\
&=&\#\{\mbox{irred.\ comps.\ of $E_K$}\}
- \#\{\mbox{irred.\ comps.\ of $E'_K$}\}\\
&=& 0.
\end{array}
$$
Applying this to $u\in R(U)$ we obtain that
the Picard number of $X'_u$ is also 1 for general $u\in R(U)$,
as required.\qed

\begin{say}[Proof of Theorem \ref{thm}]\label{pf.of.thm}
For any projective variety $X$,  $\aut(X)$ is a scheme
with countably many irreducible components
and the identity component $\aut^0(X)$ 
is a (finite dimensional) algebraic group.
For Fano varieties,
$\aut(X)$ respects the anticanonical polarization, hence
it acts faithfully on some projective embedding.
In particular,  $\aut(X)$ is a linear algebraic group.
Thus, if $\dim\aut(X)>0$ then $X$ is ruled 
(cf. \cite{ros} or \cite[IV.1.17.5]{rc-book}). 
 A ruled variety is not birationally rigid
(except for $X=\p^1$). In particular, a birationally rigid
variety of dimension $\geq 2$ has a finite automorphism group.
$\p^1$ does not satisfy any of the conditions of
(\ref{thm}.1--4). Hence from now on we may assume that
$\aut(X)$ is finite.

 Over an uncountable algebraically closed field,
a scheme which is a countable union of finite type 
subschemes has countably many points iff
it is zero dimensional. Thus 
 (1) and (4) are equivalent.

Let us prove next that (1) and (2) are equivalent.
Let $U\subset \bir(X)$ be a positive dimensional irreducible component
with generic point $z_g$. Then $z_g$ corresponds to a
birational self-map $\phi_Z$ of $X$ defined over $k(Z)$
which is not in $\bir(X)$. Conversely, if
$\phi_K\in \bir(X_K)$ is not in $\bir(X)$ then
the closure of the corresponding point in
$\hilb(X\times X)$ gives a positive dimensional component
of $\bir(X)$.

It is clear that (3) $\Rightarrow$ (2).

Finally let us prove that (1) implies (3),

Let $X'_K$ be a Mori fiber space birational to $X_K$.
We may assume that everything is defined over 
a $k$-variety $U$.
Thus we have
$$
X_U\stackrel{\phi_U}{\map} X'_U \stackrel{\pi_U}{\to} S_U.
$$
If $S_U\to U$ has positive dimensional fibers,
then for almost all $u\in U(k)$, the fiber $X_u$ is
birational to a nontrivial Fano fiber space; a contradiction.
Thus we may assume that $S_U=U$.

Here we are not allowed to replace $U$ with an arbitrary
ramified cover, hence
we can not
assert  that $X'_u$ has Picard number 1 for general $u\in U(k)$.

There is, however, a quasi-finite Galois cover $W\to U$
 such that the whole Picard group of the
geometric generic fiber of $X'_W$ is defined over $k(W)$.
Let us run the MMP on
$X'_W\to W$ to end up with
$\tau: X'_W\map X^*_W$ and a Mori fiber space
$ X^*_W\to S^*_W$.
(This much of the MMP is known in any dimension by \cite{bchm}.)
As before we see that $S^*_W\cong W$ and
the fibers of $X^*_W\to W$ have Picard number 1.
The fibers of $X^*_W\to W$ over $k$-points are thus isomorphic to $X_k$
and so, arguing with 
$\isom\bigl(X_k\times W, X^*_W\bigr)$ as above,
we obtain that 
 $X^*_W\to W$ is an \'etale locally trivial $X$-bundle
over a dense open subset of $W$.
After a possible further covering of $W$, we may assume that 
 $X^*_W\cong X\times W$.  

Note that the Galois group $G=\gal(W/U)$ need not act trivially on
the Picard group, thus we use the ordinary MMP, not the
$G$-equivariant MMP. In particular, the $G$-action on
$X'_W$ does not give a regular $G$-action on $X^*_W$,
only a birational $G$-action.

Thus, in (\ref{pf.of.thm}.1) the maps are $G$-equivariant,
the $G$-action is trivial on the $X$-factor on the left
but on  $X^*_W\cong X\times W$ we have a so far unknown
birational $G$-action which commutes with projection to $W$. 
$$
X\times W=X_W\stackrel{\phi}{\map} X'_W\stackrel{\tau}{\map} X^*_W
\cong X\times W.
\eqno{(\ref{pf.of.thm}.1)}
$$
A birational map 
$X\times W \map X\times W$ which commutes with projection to $W$
corresponds to a rational map $W\map\bir(X)$. Since
$\bir(X)$ is zero dimensional, any such 
birational map is obtained by $\bir(X)$ acting on the $X$-factor.
Thus the $G$-action on  $X^*_W\cong X\times W$
is given by 
$$
(x,w)\mapsto \bigl(\phi(g)(x), g(w)\bigr)
$$
where $\phi:G\to \bir(X)$ is a homomorphism.

Applying this argument to $\tau\circ \phi$, we get that
$$
\bigl(\tau\circ \phi\bigr)(x,w)=\bigl(\rho(x), w\bigr)
$$
for some $\rho\in \bir(X)$.

The $G$-equivariance  of $\tau\circ \phi$ gives
the equality
$$
\bigl(\rho(x), g(w)\bigr)=\bigl(\phi(g)\rho(x), g(w)\bigr)
\quad \forall g\in G.
$$
Thus $\phi(g)$ is the identity and so the 
$G$-action on  $X^*_W\cong X\times W$ is in fact 
trivial on the $X$-factor, hence biregular. 

We can now take the quotient by $G$ to obtain
$$
\tau/G:X'_U=X'_W/G\map X^*_W/G\cong X\times U.
$$
Note that $\tau$ does not extract any divisor,
thus $\tau/G$ also does not extract any divisor.
Since $X'_U$ has relative Picard number 1 over $U$,
we conclude that $\tau/G$ also does not contract
any divisor. By a lemma of Matsusaka and Mumford
(cf.\ \cite[5.6]{ksc}), this implies that
$\tau/G$ is an isomorphism.
\qed
\end{say}

 \begin{ack}  I thank I.\ Cheltsov
for useful e-mails and comments. 
Partial financial support  was provided by  the NSF under grant number 
DMS-0500198. 
\end{ack}

\providecommand{\bysame}{\leavevmode\hbox to3em{\hrulefill}\thinspace}
\providecommand{\MR}{\relax\ifhmode\unskip\space\fi MR }
\providecommand{\MRhref}[2]{%
  \href{http://www.ams.org/mathscinet-getitem?mr=#1}{#2}
}
\providecommand{\href}[2]{#2}

\bibliography{refs}

\vskip1cm

\noindent Princeton University, Princeton NJ 08544-1000

\begin{verbatim}kollar@math.princeton.edu\end{verbatim}

\end{document}